\documentclass[12pt]{amsart}

\usepackage{amsmath}
\usepackage{amsfonts, amssymb, amsthm}
\usepackage[all]{xy}

\newtheorem{theorem}{Theorem}
\newtheorem{proposition}{Proposition}
\newtheorem{lemma}{Lemma}

\theoremstyle{remark}

\newcommand{\be}{\begin{equation}}
\newcommand{\ee}{\end{equation}}
\newcommand{\bes}{\begin{equation*}}
\newcommand{\ees}{\end{equation*}}
\newcommand{\np}{\noindent}
\newcommand{\tei}{{Teichm\"uller}}

\newcommand{\RR}{\mathbb R}
\newcommand{\CC}{\mathbb C}
\newcommand{\PP}{\mathbb P}

\def\ov{\overline}
\def\del{\partial}
\def\delbar{\overline{\partial}}
\def\Box{ { \sqcup\llap{ $\sqcap$} } }

\begin{document}
\title{Teichm\"uller space via Kuranishi families}
\author{Enrico Arbarello}
\address{Dipartimento di Matematica ``G. Castelnuovo'', Universit\`a di Roma ``La Sapienza'', p.le A. Moro 2, 00185 Roma, Italy}
\email{ea@mat.uniroma1.it}
\thanks{Research partially supported by: PRIN 2005 \textit{Spazi di moduli e teoria di Lie}; FAR 2004 (Pavia) \textit{Variet\`a algebriche, calcolo algebrico, grafi orientati e topologici}}

\author{Maurizio Cornalba}
\address{Dipartimento di Matematica, Universit\`a di Pavia, Via Ferrata 1, 27100 Pavia, Italy}
\email{maurizio.cornalba@unipv.it}

\subjclass[2000]{30F60, 14H15 (primary); 32G15, 14H10 (secondary).}

\begin{abstract}
In this partly expository note we construct {\tei} space by patching together Kuranishi families. We also discuss the basic properties of {\tei} space, and in particular show that our construction leads to simplifications in the proof of {\tei}'s theorem asserting that the genus $g$ {\tei} space is homeomorphic to a $(6g-6)$-dimensional ball.
\end{abstract}

\maketitle
\section{Introduction}
Our main goal, in this partly expository note, is to show how the genus $g$ {\tei} space $T_g$ and the universal family over it can be constructed by patching together Kuranishi families of genus $g$ curves. This approach, which is close in spirit to the one of Grothendieck \cite{Grothendieck_1,Grothendieck_2}, should appeal in particular to readers with an algebro-geometric background, as it relies mostly on standard tools and methods of their trade. In fact, the main objects we use are the Kuranishi families of smooth genus $g$ curves, whose construction and properties can be derived quite directly from the theory of the Hilbert scheme. Our point of view seems to have several advantages over more traditional ones. First of all, {\tei} space is constructed directly as a complex manifold, and the construction makes it obvious that it enjoys a natural universal property. Secondly, it is very easy to show that the action of the {\tei} modular group on $T_g$ is properly discontinuous. Finally, our presentation provides a shortcut to the proof of Teichm\"uller's theorem, which states that $T_g$
is homeomorphic to a ($6g-6$)-dimensional ball $B_g$. In fact, we show that the universal property of the Kuranishi family holds not only with respect to analytic deformations, but also with respect to {\it continuous} ones, and observe that this immediately proves that the {\tei} map $\Phi\colon B_g\to T_g$ is continuous. This first step eliminates some of the technicalities involved in the usual proof of {\tei}'s theorem, highlighting the two central points of the argument: {\tei}'s uniqueness theorem and the existence of solutions of the Beltrami equation. It should be observed that the same program can be carried out for the {\tei} space $T_{g,n}$ of genus $g$ curves with $n$ marked points, although this will not be done here. In the bibliography we have provided references to the most fundamental texts in this very classical theory.

\section{The Kuranishi family}\label{kuranishi}
Throughout this section, $C$ will denote a compact connected Riemann surface of genus $g>1$, or more briefly, a {\it curve of genus $g>1$.} A {\it family of genus $g$ curves parametrized by an analytic space $S$} is a proper surjective analytic map
$\eta\colon X\to S$ having genus $g$ compact Riemann surfaces as fibers.
If $s$ is a point of $S$, we shall write $X_s$ to indicate the fiber $\eta^{-1}(s)$; more generally, we shall use the convention of appending the subscript $s$ to an object over $S$ to indicate its ``fiber'' at $s$. The precise meaning of this in each case will be clear from the context.
We will say that $\eta\colon X\to S$ is a {\it family of genus $g$ curves in $\mathbb P^N$, parametrized by $S$}, if
$ X$ is a subvariety of $\mathbb P^N\times S$ and $\eta$ is the restriction to $ X$ of the projection on the second factor.
A {\it deformation of $C$, parametrized
by a pointed analytic space $(S, s_0)$,} is a family of genus $g$ curves
$\eta\colon X\to S$, together with an
isomorphism $\psi\colon C\overset\simeq\to X_{s_0}$.
Thus, a deformation of $C$ is determined by the data $(\eta, \psi)$.
Given a neighborhood $U$ of $s_0$, we let $\eta_U\colon \eta^{-1}(U)\to U$ be the map obtained by restriction from $\eta$.
If $\eta'\colon {X}'\to S',\ \psi'\colon C\to X'_{s'_0}$ is another deformation of $C$, parametrized by a pointed analytic space $(S',s'_0)$, a {\it morphism of deformations} from $(\eta,\psi)$ to $(\eta',\psi')$ is a pair $(f,F)$ of morphism fitting in a cartesian diagram
$$
\xymatrix{ 
X \ar[d]_\eta \ar[r]^F 
&X' \ar[d]_{\eta'} \\ 
(S,s_0) \ar[r]^f & (S',s_0') }
$$
Let $(B,b_0)$ be a $(3g-3)$-dimensional connected pointed complex manifold.
A {\it Kuranishi family for $C$ } is a deformation of $C$ parametrized by $(B,b_0)$:
\be
\pi\colon \mathcal C\to B\,,\quad \varphi\colon C\overset\simeq\to\mathcal C_{b_0}\label{kf}
\ee
satisfying the following universal property.
Given a deformation $(\eta, \psi)$ of $C$ as above, there are a
neighborhood $U$ of $s_0$ and a morphism $(f,F)$ of deformations
from $(\eta_U, \psi)$ to $(\pi,\varphi)$; moreover, this morphism is essentially unique, in the sense that any other morphism of this sort agrees with $(f,F)$ over $\eta^{-1}(U')\to U'$, for some neighborhood $U'$ of $s_0$.

Kuranishi families for curves of genus $g>1$ always exist. We briefly illustrate how one such family may be constructed.
Let $\nu\geq 3$ and set
$N=(2\nu-1)(g-1)-1$.
As is well known, a curve of genus $g$ can be embedded in the projective 
space $\mathbb P^{N}$ via the $\nu$-canonical map, that is, via the sections of the $\nu$-th power of the canonical line bundle $\omega_C$. The embedding depends on the choice of a basis for the space of sections of $\omega_C^\nu$, and different choices of basis give rise to projectively equivalent embeddings. Conversely, any isomorphism between two genus $g$ curves comes from a unique projectivity of $\mathbb P^N$ carrying the $\nu$-canonical image of one curve to the $\nu$-canonical image of the other.
The {\it Hilbert scheme $H_{g,\nu}$ of
$\nu$-canonical, genus $g$ curves } is a smooth projective variety
whose points are in a one-to-one correspondence with the set
of $\nu$-canonically embedded curves of genus $g$. The Hilbert scheme
$H_{g,\nu}$ is the parameter space for a family $\xi\colon \mathcal Y\to H_{g,\nu}$ of genus $g$ curves in $\mathbb P^{N}$
having the following universal property.
For every family $\eta\colon X\to S$ of genus $g$ curves $\nu$-canonically embedded in
$\mathbb P^{N}$ there exists a unique morphism $\alpha\colon S\to H_{g,\nu}$
exhibiting the family $\eta$ as the pull-back, via $\alpha$,
of the universal family $\xi$, so that there is a commutative diagram of cartesian squares
$$
\xymatrix@H+15pt{ 
\mathbb P^N\times S \ar[r]^{{\bf 1}\times \alpha}
&\mathbb P^N\times H_{g,\nu}\\ 
X \ar[r]^h \ar@{^{(}->}[u] \ar[d]_\eta & \mathcal Y \ar @{^{(}->}[u] \ar[d]^\xi \\
S \ar[r]^\alpha &H_{g,\nu} }
$$
The projective group $PGL(N+1,\CC)$ acts naturally on
$H_{g,\nu}$, with finite stabilizers, since curves of genus $g>1$ have finite automorphism groups. Thus, in accordance with Riemann's count,
$$
\dim H_{g,\nu}=3g-3 +\dim PGL(N+1,\CC)\,.
$$
Now let us consider a genus $g$ curve $C$ and a $\nu$-canonical embedding
$\varphi\colon C\overset\simeq\to\Gamma\subset\mathbb P^N$ of $C$, associated to a basis
$\chi_0,\dots,\chi_N$ of the space of $\nu$-fold holomorphic differentials on $C$. Let $p\in H_{g,\nu}$ be the point corresponding to
$\Gamma$. We choose
a local ``slice'' at $p$ for the action of $PGL(N+1,\CC)$, that is, a
$(3g-3 )$-dimensional locally closed subvariety $U$ passing through $p$ and transverse to the orbits of $PGL(N+1,\CC)$.
Possibly after shrinking $U$, if $W$ is a sufficiently small open neighborhood of the identity in the projective group, $(g,u)\mapsto gu$ is an isomorphism between $W\times U$ and an open neighborhood $V$ of $p$. Thus there is a canonical fibration $\sigma\colon V\to U$ along the orbits of $PGL(N+1,\CC)$.
Let now $\pi\colon \mathcal C\to U$ be the restriction to $U$ of the universal family $\xi\colon \mathcal Y\to H_{g,\nu}$. This
family, together with the identification $\varphi\colon C\to\Gamma=\pi^{-1}(p)$,
is a Kuranishi family for $C$. Let us see why. Suppose
$$
\eta\colon X\to S\,,\qquad\psi\colon C\to X_{s_0}
$$
is a deformation of $C$ parametrized by $(S,s_0)$. Consider the sheaf $\omega_\eta$ of relative holomorphic differentials along the fibers of $\eta$.
By suitably shrinking the parameter space
$S$ of this deformation we may assume that $\eta_*\omega_\eta^\nu$
is trivial on $S$. Choosing a frame of this vector bundle
exhibits $\eta\colon X\to (S,s_0)$ as a family of $\nu$-canonically embedded curves, that is, determines an embedding $\beta\colon X\to \mathbb P^N\times S$ whose fiber $\beta_s\colon X_s\to\mathbb P^N$ at $s$ is $\nu$-canonical for every $s\in S$.
Furthermore, if the frame is chosen so as to pull back to $\chi_0,\dots,\chi_N$ on $C$, then the composition of the identification $\psi\colon C\to X_{s_0}$ with $\beta_{s_0}\colon X_{s_0}\to \mathbb P^N$ is equal to the composition of $\varphi$ with the inclusion of $\Gamma$ in $\mathbb P^N$. The universal property of the Hilbert scheme then yields a morphism $\alpha\colon S\to H_{g,\nu}$ such that $\alpha(s_0)=p$.
We may assume that $\alpha(S)\subset V$. It is then immediate to verify
that the morphism $\sigma\alpha\colon S\to U$ realizes the universal property
of a Kuranishi family.

Possibly shrinking $B$, one can assume that the Kuranishi family (\ref{kf}) satisfies various useful additional properties. For instance, one can set things up so that (\ref{kf}) is ``Kuranishi at every point of $B$''. This means that, for any $b\in B$, the deformation consisting of $\mathcal C\to B$ and of the identity isomorphism from $\mathcal C_b$ to itself is a Kuranishi family for $\mathcal C_b$. In fact, the Kuranishi family whose construction we outlined above has this property.

We may also assume that the automorphism group of $C$ acts on $\mathcal C$ and $B$. In fact, let $\gamma$ be an automorphism of $C$. Then
\be
\pi\colon \mathcal C\to B\,,\quad \varphi\gamma\colon C\overset\simeq\to \mathcal C_{b_0}\label{kf2}
\ee
is another Kuranishi family. By the universal property, up to shrinking $B$, if necessary, there is a unique cartesian diagram
$$
\xymatrix{
\mathcal C \ar[r]^{\widetilde{f}_\gamma}
\ar[d]_\pi & \mathcal C \ar[d]^\pi \\
(B,b_0) \ar[r]^{f_\gamma}&(B,b_0)\\
}
$$
where $f_\gamma$ and $\widetilde{f}_\gamma$ are isomorphisms,
inducing an isomorphism between the two deformations of $C$ given by
(\ref{kf}) and (\ref{kf2}).
This results, possibly after further shrinking of $B$, in compatible actions of $\operatorname{Aut}(C)$ on $B$ and $\mathcal C$ pulling back, via $\varphi$, to the standard action on $C$. The action on $\mathcal C$ is always faithful, while the action on $B$ is not faithful only when $g=2$; in this case the only elements of $\operatorname{Aut}(C)$ acting trivially are the identity and the hyperelliptic involution.
By further shrinking the base $B$ one may even assume that any isomorphism between two fibers
$\mathcal C_b$ and $\mathcal C_{b'}$ of $\pi$ comes by restriction from $\widetilde{f}_\gamma$, for some $\gamma\in \operatorname{Aut}(C)$. This is often crucial in applications.

We end this section by recalling a fundamental compactness property of families of compact Riemann surfaces. Suppose $\alpha\colon X\to U$ and $\alpha'\colon X'\to U'$ are two families of compact Riemann surfaces of genus $g>1$ parametrized by analytic spaces $U$ and $U'$. Let $\{u_n\}_{n\in \mathbb N}$ ($\{u'_n\}_{n\in \mathbb N}$) be a sequence of points in $U$ (resp., $U'$) converging to a point $u\in U$ (resp., $u'\in U'$). Assume that $X_{u_n}$
is isomorphic to $X'_{u'_n}$, for every $n$. Then $X_{u}$ is isomorphic to $X'_{u'}$.
This can be deduced from the following purely algebro-geometric statement, which is in turn a rather straightforward consequence of the properness of relative Hilbert schemes and of the theory of minimal models of algebraic surfaces. Let $Y\to R$ and $Y'\to R$ be {\em algebraic} families of smooth curves of genus $g>1$; then the scheme $\operatorname{Isom}_R(Y,Y')$ parametrizing isomorphisms $Y_r\overset\simeq\to Y'_r$, $r\in R$, is proper over $R$.
Here is how the reduction to this statement of the preceding one goes. Let $Z\to S$ and $W\to T$ be algebraic families of smooth curves of genus $g>1$. Set $Y=Z\times T$, $Y'=S\times W$, $R=S\times T$. Then the locus $\{(s,t):X_s\text{ is isomorphic to } Y_t\}\subset S\times T=R$ is the projection of $\operatorname{Isom}_R(Y,Y')$, and hence is closed. This shows in particular that the statement to be proved is true for algebraic families. To prove it in general it then suffices to show that the original families $\alpha\colon X\to U$ and $\alpha'\colon X'\to U'$, or just their restrictions to suitable neighborhoods of $u$ and $u'$, come by pullback from algebraic families. But this is clear. In fact, up to shrinking $U$ and $U'$, we may suppose that $\alpha_*\omega^\nu_\alpha$ and $\alpha'_*\omega^\nu_{\alpha'}$ are trivial for some $\nu\ge 3$. Thus we may embed $X$ and $X'$ $\nu$-canonically, so $X\subset \mathbb P^N\times U$, $X'\subset \mathbb P^N\times U'$, and both families come from the universal family over the Hilbert scheme $H_{g,\nu}$ via maps $U\to H_{g,\nu}$ and $U'\to H_{g,\nu}$.

A convenient reference for deformation theory and Hilbert schemes is \cite{Sernesi}. A detailed presentation of Kuranishi families will be contained in the forthcoming volume \cite{ACGHII}.

\section{The construction of Teichm\"uller space as a complex manifold}\label{teichconstr}

Fix an oriented genus $g$ topological surface $\Sigma$. We assume that
$g>1$. Given a genus $g$ curve $C$, a {\it {\tei} structure}
on it is the isotopy class $[f]$ of an
orientation-preserving homeomorphism
$
f\colon C\to \Sigma\,.
$
An isomorphism between curves with {\tei} structure $(C, [f])$ and
$(C',[f'])$ is an isomorphism of curves
$
\varphi \colon C\to C'\,,
$
such that $[f'\varphi]=[f]$. The set of isomorphism classes of
genus $g$ curves with {\tei} structure has a natural topology and
complex structure which we are presently going to describe. The resulting space is
called the {\it {\tei} space of $\Sigma$} and is denoted by the symbol
$T_\Sigma$. The point in $T_\Sigma$ associated to the curve $C$
and to the isotopy class $[f]$ will be denoted by the symbol
$[C,[f]]$. If $\Sigma'$ is another oriented
genus $g$ surface, and $\rho\colon \Sigma\to \Sigma'$ an oriented homeomorphism, we get a bijection $T_\Sigma\to T_{\Sigma'}$
by sending $[C, [f]]$ to $[C, [\rho\circ f]]$. Clearly, this bijection depends only on the isotopy class of $\rho$, and it will be clear from the construction of the complex structure on {\tei} space that it is an isomorphism of complex manifolds. We are therefore justified in writing $T_g$, instead of $T_\Sigma$, when the reference surface $\Sigma$ is kept fixed, as will usually happen.

It is important to remark that a curve with {\tei} structure $(C,[f])$
has no non-trivial automorphism.
In fact, an automorphism $\varphi\colon C\to C$
such that
$[f\varphi]=[f]$ must be homotopically trivial and, in particular,
must induce the
identity on complex cohomology. This implies that $\varphi$ induces the identity at the level of holomorphic forms, and hence must commute with the canonical map. This shows that $f$ is the identity if $C$ is not hyperelliptic, and that $f$ is the identity or the hyperelliptic involution otherwise. However, the hyperelliptic involution does not act trivially on complex cohomology, so $f$ must be the identity in this case as well.

We next introduce a topology and a complex structure on {\tei} space.
We begin by extending the
notion of {\tei} structure from curves to families of curves. A {\it {\tei} structure} on a family
of
genus
$g$ curves $\eta\colon X\to S$ is the datum of a
{\tei} structure $[f_s]$ on each fiber $X_s$, satisfying the following coherence condition. There exists a cover of $S$
with open sets $U$ together with topological trivializations
$(F,\eta_U)\colon \eta^{-1}(U)\overset\simeq\to \Sigma\times U$ such that, for each $s\in U$, $[F_s]=[f_s]$. Let $V$ be an open subset of $S$, and let
$$
\aligned
(F,\eta_V)\colon &\eta^{-1}(V)\overset\simeq\to \Sigma\times V\\
(G,\eta_V)\colon &\eta^{-1}(V)\overset\simeq\to \Sigma\times V
\endaligned
$$
be two topological trivializations. Suppose that $[F_s]=[G_s]$ for some $s\in V$. Let $\gamma\colon [0,1]\to V$ be a path from a point $s'\in V$ to $s$.
Then, as $t$ varies between 0 and 1, the homeomorphisms
$$
G_{\gamma(t)}\circ F_{\gamma(t)}^{-1}\circ
F_{s'}\colon X_{s'}\to \Sigma
$$
describe an isotopy between $G_{s'}$ and $G_s\circ F_s^{-1}\circ F_{s'}$,
which in turn is isotopic to $F_{s'}$. Hence $[F_{s'}]=[G_{s'}]$ for all $s'\in V$ belonging to the same connected component as $s$.
This implies in particular that, when $V$ is connected, $F$ and $G$ are equal if and only if they agree at one point.

Now let $(C,[f])$ be a genus $g$ curve with {\tei} structure.
Consider a Kuranishi family (\ref{kf}) for the curve $C$.
Possibly after shrinking $B$, such a family admits a topological trivialization
$(F,\pi)\colon \mathcal C\to \Sigma\times B$ such that $F_{b_0}\circ \varphi=f$, and hence can
be endowed with a unique {\tei} structure extending the one on
$C$. We thus get what we shall refer to as a {\it Kuranishi
family for the curve with {\tei} structure} $(C, [f])$. The name is
justified by the fact that such a family enjoys, with respect to deformations of
curves with {\tei} structure, a universal property exactly analogous to the one of
standard Kuranishi families, as follows immediately from the universal
property of ordinary Kuranishi families and the uniqueness of the {\tei} structure
extending the one on the fiber at $s_0$. It should be observed that, when $B$ is
small enough, the family we just constructed is Kuranishi at every point of $B$, as
follows from the analogous property of standard Kuranishi families.

We are ready to describe the topology and the complex
structure on $T_g$. Let $y=[C,[f]]$ be a point of $T_g$. Choose a Kuranishi family for
$(C, [f])$. As we explained, this can be constructed by putting on
the family (\ref{kf}), where $B$ is chosen to be connected and
``sufficiently small'', the unique {\tei} structure extending $[f]$; for each
$b\in B$, we denote by $[F_b]$ the {\tei} structure on
$C_b$. As we pointed out at the end of Section \ref{kuranishi}, we may
suppose that the action of $G=\operatorname{Aut}(C)$ on
$C$ extends to an equivariant action on $\mathcal C$ and on $B$, and that any isomorphisms between fibers of $\pi$ is
the restriction of the action of an element of $G$ on $\mathcal C$. We claim that the natural map
$$
\aligned
\alpha \colon &\,B \longrightarrow T_g\\
&\,b\mapsto [C_b, [F_b]]
\endaligned
$$
is injective. A map of this sort will be called a {\it standard coordinate
patch} for $T_g$ around the point $y$.
To prove the injectivity of $\alpha$, suppose that $\alpha(b)=\alpha(b')$ for $b,b'\in B$. Thus there is an isomorphism
$\psi\colon C_b
\to C_{b'}$ such that
$[F_b]=[F_{b'}\circ\psi]$. Our assumptions imply that $\psi$ is induced by an
element $\rho\in G$, so that there is a commutative diagram
\be
\xymatrix{
\mathcal C \ar[r]^{\tilde\xi} \ar[d]_\pi
&\mathcal C \ar[d]^\pi \\
B \ar[r]^{\xi}& B
}
\label{autaction}
\ee
where $\tilde\xi$ and $\xi$ are automorphisms such
that $\xi(b_0)=b_0$, $\xi(b)=b'$, $\tilde\xi_b=\psi$, and $\tilde\xi$ pulls back to $\rho$ on $C$. We define a new {\tei} structure on (\ref{kf}) by setting $F'_t=F_{\xi(t)}\circ\tilde\xi_t$. Since
$[F'_b]=[F_{b'}\circ\psi]=[F_b]$, it follows that $[F'_t]=[F_t]$ for every $t\in
B$. When $t=b_0$, this says that $[f\circ\rho]=[f]$. Thus
$\rho$ is isotopic to the identity and, as we observed above, this implies that
$\rho=1$, proving our claim.

Let us show that the patches we have just introduced define a complex structure on
$T_{g}$. Let
$\alpha\colon B\to T_g$ and $\beta\colon B'\to T_{g} $ be two standard patches whose
codomains have a point $z$ in common. Let $\pi\colon \mathcal C\to B$, $\pi'\colon \mathcal C'\to B'$ be the corresponding families, and $U$ a small enough neighborhood of $\beta^{-1}(z)$. Then, by the universal property of Kuranishi families, there is a unique morphism of families of curves with {\tei} structure
$$
\xymatrix{
{\pi'}^{-1}(U) \ar[r]^{\Psi} \ar[d]_{\pi_U'}
& \mathcal C \ar[d]^\pi\\
U \ar[r]^{\psi} & B
}
$$
On the other hand, $\psi$ clearly agrees with the restriction to $U$ of
$\alpha^{-1}\circ\beta$. In particular, $\beta^{-1}(\alpha(B))$ is open in $B'$, and $\alpha^{-1}\circ\beta$ is holomorphic on it.

This completes the construction of a (possibly non-Hausdorff) complex structure on $T_g$; in
the next section we shall prove that $T_g$ is actually a Hausdorff topological space. In fact, the above argument also shows that the Kuranishi families corresponding to standard coordinate patches can be canonically glued together to yield a {\it universal
family of curves with {\tei} structure}
\bes
\eta\colon \mathcal X_{g} \longrightarrow T_g\,.
\ees
This family is universal in the sense that any holomorphic family of genus $g$ curves with
{\tei} structure over a base $T$ is isomorphic, via a unique isomorphism,
to the pull-back of the family $\mathcal X_{g}\to T_g$ via a unique morphism
$T\to T_g$. In particular, $T_g$ represents the functor
$$
T\mapsto\left\{\genfrac{}{}{0pt}{}{\text{isomorphism classes of holomorphic families of}}{\text{genus $g$ curves with {\tei} structure over $T$}}\right\}
$$
from analytic spaces to sets.

\section{The mapping class group and its action}\label{modular}

As in the previous section, we pick a reference oriented genus $g$ topological surface $\Sigma$. The {\it mapping class group} $\Gamma_{\Sigma}$, also called
{\it {\tei} modular group}, is the group of all isotopy classes of orientation-preserving
homeomorphism of $\Sigma$. When the reference surface $\Sigma$ is kept fixed,
we shall often denote it by $\Gamma_g$.
The mapping class group acts naturally on
$T_g= T_{\Sigma}$, the action of an element
$[\gamma]$ being given by
$$
[\gamma]\cdot[C, [f]]= [C, [\gamma\circ f]]\,.
$$
It is immediate to show that the elements of $\Gamma_g$ act on $T_g$ as holomorphic automorphisms. Consider an element $[\gamma]\in \Gamma_g$, and a standard coordinate patch $\alpha\colon B\to T_g$ obtained from a Kuranishi
family (\ref{kf}), endowed with a {\tei} structure $\{F_b\}_{b\in B}$. We know that $[\gamma]$ acts by replacing each
{\tei} structure $[f]$ with $[\gamma\circ f]$.
Substituting $\{\gamma\circ F_b\}_{b\in B}$ for the {\tei} structure $\{F_b\}_{b\in B}$ produces a new
coordinate patch $\beta\colon B\to T_g$, and it is clear that, for any $b\in B$,
$$
\beta^{-1}([\gamma]\cdot \alpha(b))=b\,.
$$
This shows that $[\gamma]$ acts holomorphically.

As a set, the \it moduli space $M_g$ of genus $g$ curves \rm is the set of isomorphism
classes of smooth curves of genus $g$. Ignoring {\tei} structures
gives a map
$$
m\colon T_g\to M_g\,,
$$
which can be identified with the quotient map from
$T_g$ to $T_g/\Gamma_g$. To see this first observe that, given $y\in T_g$ and
$[\gamma]\in\Gamma_g$, $y$ and $[\gamma]\cdot y$ obviously map to the same point of $M_g$. Conversely, given points $y=[C, [f]]$
and $y'=[C', [f']]$ of $T_g$ mapping to the same point of $M_g$, there is an isomorphism $\varphi\colon C\to C' $, so that $y'=[\gamma]\cdot y$, where $\gamma=f'\varphi f^{-1}$. Let now $G$ be the automorphism group of the curve $C$. We claim that the homomorphism $\sigma$ from $G$ to $\Gamma_g$ given by
$$
\rho\mapsto[\gamma_\rho]\,\text{, where\ \ }\gamma_\rho=f\rho f^{-1}\,,
$$
identifies $G$ with the stabilizer of $y$. First of all, $\sigma$ is injective, since $[f\rho f^{-1}]=1$ if and only if $\rho$ is homotopic to the identity, and, as we know, this happens only when $\rho=1$. Furthermore, $f\rho=\gamma_\rho f$, which exactly means that $\rho$ is an isomorphism
between the curves with {\tei} structure $(C, [f])$ and
$(C, [\gamma_\rho f])$. Thus $y=[\gamma_\rho]\cdot y$, that is, $\sigma(\rho)$ belongs to the stabilizer of $y$. Finally, given any $[\gamma]$ in the stabilizer of $y$, there is an element $\psi$ of $G$ such
that $[f\psi]=[\gamma f]$, that is, $[\gamma]=[\gamma_\psi]=\sigma(\psi)$ for some $\psi\in G$.

Now consider a standard coordinate
patch $\alpha\colon B\to T_g$ centered at $y$, coming from a Kuranishi family
(\ref{kf}), endowed with a {\tei} structure $\{F_b\}_{b\in B}$. We shall prove that $\alpha$ is $G$-equivariant. If $\rho$ is
an element of $G$, there is a diagram (\ref{autaction}), where $\tilde\xi$
and $\xi$ are automorphisms such that $\xi(b_0)=b_0$ and
$\tilde\xi_{b_0}\varphi=\varphi\rho$.
The action of $\rho$ on $B$ is then
given by $\rho\cdot b=\xi(b)$. We must show that
$\alpha(\xi(b))=[\gamma]\cdot \alpha(b)$, where
$\gamma=\gamma_\rho=f\rho f^{-1}$. We define two new {\tei} structures
$\{H_b\}_{b\in B}$ and $\{K_b\}_{b\in B}$ on (\ref{kf}) by setting
$H_b=F_{\xi(b)}\tilde\xi_b$ and $K_b=\gamma F_b$. It is easy to verify that $H_{b_0}=K_{b_0}$. It follows that $[H_b]=[K_b]$ for all
$b$. Since $\tilde\xi_b$ is an isomorphism between $(C_b, [H_b])$ and
$(C_{\xi(b)}, [F_{\xi(b)}])$, and since
$$
[\gamma]\cdot[C_b, [F_b]]
=[C_b, [K_b]]\,,
$$
we have that
$$
[\gamma_\rho]\cdot\alpha(b)=[C_b, [K_b]]
=[C_b, [H_b]]=\alpha(\xi(b))\,,
$$
as wished.

The moduli space $M_g=T_g/\Gamma_g$ is equipped with the quotient topology. A basis for this topology is given as follows.
As we just saw, any standard patch $\alpha\colon B\to T_g$ drops to an injective map $\overline \alpha\colon B/G\to M_g$, where $G$ is the automorphism group of the central fiber. The open sets of the form
$\overline \alpha(B/G)$ form a basis for the topology of $M_g$.

The fact that the moduli space $M_g$ is Hausdorff is an immediate consequence of the property mentioned at the end of Section \ref{kuranishi}.
Suppose in fact that the sequence $\{z_n\}\subset M_g$ converges to two points $x$ and $y$, corresponding to the isomorphism classes
of curves $C$ and $C'$. Consider Kuranishi families
$\pi\colon \mathcal C\to B$
and $\pi'\colon \mathcal C'\to B'$ for $C=\pi^{-1}(b_0)$ and $C'=\pi^{-1}(b'_0)$, respectively, yielding standard patches $\alpha \colon B\to T_g$ and $\alpha'\colon B'\to T_g$. We may assume that
$\{z_n\}\subset \overline{\alpha}(B)\cap\overline{\alpha}'(B')$.
Lift the sequence
$\{z_n\}$, via $\overline{\alpha}$ and $\overline{\alpha}'$, to sequences $\{x_n\}\subset B$ and $\{y_n\}\subset B'$ converging, respectively, to $b_0$ and $b'_0$. Since $C_{x_n}$ is isomorphic to $C'_{y_n}$, for all $n$, we conclude that $C\cong C_{b_0}$ is isomorphic to $C'\cong C'_{b'_0}$, so that $x=y$.

We shall now show that $T_g$ is Hausdorff and, at the same time,
that the action of $\Gamma_g$ on $T_g$ is properly discontinuous.
We must prove that:
\begin{itemize}
\item[ 1)] Points $y$ and $y'$ in $T_g$ belonging
to different $\Gamma_g $-orbits possess disjoint neighborhoods.
\item[2)] Every point $y\in T_g$ possesses a neighborhood $V$ such that,
if $G_y$ is the stabilizer of $y$, then
\bes
\{\gamma\in \Gamma_g:\gamma V\cap V\neq \emptyset\}\subset G_y\,.
\ees
\end{itemize}
\np
The first property is a direct consequence of the Hausdorffness of $M_g$.
To prove the second, fix
$y\in T_g$.
We can then identify the stabilizer $G_y$ with the automorphism group of the fiber $\mathcal X_y$ of the universal family $\mathcal X\to T_g$. We claim that $y$ has arbitrarily small neighborhoods $V$ with the following
properties:
\begin{itemize}
\item[\sl a)] $V$ is $G_y$-stable;
\item[\sl b)] for $y$ and $z$ in $V$, the fiber $\mathcal X_y$ is isomorphic to the fiber $\mathcal X_z$
if and only if $z=gy$, for some $g\in G_y$;
\item[\sl c)] the stabilizer of $z$ is contained in $G_y$ for any $z\in V$.
\end{itemize}

In fact, we can take as $V$ the images $\alpha(B)$ of standard patches
$\alpha\colon B\to T_g$ centered at $y$, corresponding to Kuranishi families $\mathcal C\to B$. The only one of the above properties that may not be clear is
perhaps {\sl c)}. To prove it, recall that the action of $G_y$ on $V$ corresponds
to the action on the base $B$ and total space $\mathcal C$ of the Kuranishi family of the automorphism group of the central fiber, and that every
isomorphism between fibers of $\mathcal C\to B$ is induced by an automorphism
of the central fiber. In particular, this is true of the automorphisms of the fiber
above a point $b\in B$. Since the group of these automorphisms corresponds to
the stabilizer of $z=\alpha(b)$, property {\sl c)} follows.

We claim that a neighborhood satisfying properties {\sl
a)}--{\sl c)} above also satisfies 2). Suppose in fact that $z=\gamma\cdot z'$, where $z,z'\in V$. It follows from {\sl b)} that there is $\rho\in G_y$ such that
$z'=\rho\cdot z$. But then $\gamma\rho\in G_z\subset G_y$, and hence $\gamma$ belongs to $G_y$, as desired. We may now show that $T_g$ is Hausdorff. We already know that points $y$ and $y'$ of $T_g$ belonging to different
orbits of $\Gamma_g$ can be separated, and property 2) shows how to separate points belonging to the same orbit.

This completes the proof of the fact that $T_g$ is Hausdorff and that the
mapping class group $\Gamma_g$ acts properly discontinuously on $T_g$.

\bigskip
We end this section by briefly discussing the definition of
{\tei}
space.
It is known that, if $X$ and $Y$ are surfaces without boundary, compact or
not, two homeomorphisms $f_1,f_2\colon X\to Y$ are homotopic if
and only if they are isotopic. A proof can be found in \cite{Epstein} (cf. in particular
Theorems 6.4 and A4). An immediate consequence is that a {\tei} structure on an
unpointed curve $C$ can also be defined as the {\it homotopy} class of an orientation-preserving
homeomorphism from $C$ to a reference surface $\Sigma$. Similarly,
$\Gamma_g$ can be defined as the group of orientation-preserving homeomorphisms of
$\Sigma$ modulo homotopy.

Another variant of the definitions of {\tei} space and mapping class group, entirely
equivalent to the original one, can be obtained by fixing a differentiable structure on
the reference surface $\Sigma$ and replacing the words
``homeomorphism'' and ``isotopy'' with ``diffeomorphism'' and ``differentiable
isotopy'' throughout. What this amounts to saying is that every class in
$\Gamma_g$
contains a diffeomorphism, and that two diffeomorphisms which are isotopic are also differentiably isotopic. Thus the
{\tei} space $T_g$ is the set of isomorphism classes of objects
$(C, \varphi)$, where $C$ is a curve of
genus $g$ and $\varphi$ is the differentiable isotopy class of an oriented
diffeomorphism from $C$ to $\Sigma$. Similarly,
$$
\Gamma_g=\operatorname{Diff}_+(\Sigma)
/\operatorname{Diff}_0(\Sigma)\,,
$$
where $\operatorname{Diff}_+(\Sigma)$ stands for the group of
orientation-preserving diffeomorphisms of $\Sigma$
and $\operatorname{Diff}_0(\Sigma)$ for its identity component, which is
nothing but the group of diffeomorphisms of $\Sigma$ which
are differentiably isotopic to the identity.

\section{Continuous families of Riemann surfaces}\label{continuous}

Roughly speaking, the universal property of the Kuranishi family of a smooth
curve $C$ can be expressed by saying that any small deformation of $C$ can
be obtained by pullback from the Kuranishi family. Our main goal in this
section is to show that the same is true for ``continuous deformations'' of $C$.

First, we need some preliminaries. Let $\alpha\colon X\to S$ be a surjective morphism of $C^m$ manifolds, where $m=0,\dots,\infty$. Suppose we can
find an open cover $\{U_i\}$ of $X$ and $C^m$ diffeomorphisms $\Phi_i=\Phi_{i,1}\times\Phi_{i,2}\colon U_i\to V_i\times W_i$ of manifolds over $S$, where $V_i$ is open in $\RR^n$ and $W_i$ is an open set in $S$, such that:
\begin{itemize}
\item[\sl a)]for any choice of $i$ and $j$, and for any $w\in W_j$, the function $v\mapsto \Phi_{i,1}\circ \Phi_j^{-1}(v,w)$ is of class $C^\infty$.
\item[\sl b)]the composition $\Phi_{i,1}\circ \Phi_j^{-1}$ and all its derivatives, of any order, with respect to the the $V_j$ coordinates, are of class $C^m$, for any choice of $i$ and $j$.
\end{itemize}
We shall then say that the atlas $\mathcal U=\{(U_i,\Phi_i)\}$ defines on $\alpha\colon X\to S$ a structure of $C^m$ {\em family of differentiable manifolds}. We shall refer to the components of the functions $\Phi_{i,1}$ as {\em vertical coordinates} and to derivatives with respect to them as {\em vertical derivatives}. We shall say that a function on an open set of $X$ is {\em adapted} if it is $C^\infty$ in the vertical coordinates and its vertical derivatives, of any order, are $C^m$; clearly, this notion does not depend on the chart in $\mathcal U$ with respect to which vertical derivatives are computed. Another atlas $\mathcal U'$ will be considered equivalent to $\mathcal U$ if the adapted functions with respect to it are the same as the adapted functions with respect to $\mathcal U$. This is equivalent to asking that the vertical coordinates in either atlas be adapted with respect to the other atlas. We shall say that equivalent atlases define on $\alpha\colon X\to S$ the same structure of $C^m$ family of differentiable manifolds. Given a $C^m$ family of differentiable manifolds, the charts of any atlas defining the $C^m$ family structure will be said to be adapted. A {\em morphism} from a $C^m$ family of differentiable manifolds $\alpha\colon X\to S$ to another family $\beta\colon Y\to T$ is a commutative square
$$
\xymatrix{
X \ar[r]^F \ar[d]_\alpha & Y \ar[d]^\beta \\
S \ar[r]^f & T \\
}
$$
where $f$ and $F$ are $C^m$ morphisms and the composition of $F$ with any vertical coordinate on $Y$ is an adapted function. Clearly, the class of morphisms of $C^m$ families of differentiable manifolds is closed under composition.

Let $\alpha\colon X\to S$ be a $C^m$ family of differentiable manifolds. There is an obvious notion of {\em $C^m$ family of differentiable vector bundles} (real or complex) on $\alpha\colon X\to S$. Examples are provided, for instance, by the relative tangent bundle to $X\to S$ and its dual. The class of $C^m$ families of differentiable vector bundles is closed under the standard vector bundle operations, such as passing to the dual, direct sum, tensor product, and exterior power. If $E$ is a $C^m$ family of differentiable vector bundles, it makes sense to speak of adapted local trivializations and {\em adapted sections} of $E$. In particular, it makes sense to speak of adapted relative differential forms (or, more generally, relative $E$-valued differential forms) along the fibers of $\alpha$. There is also a good notion of $C^m$ family of linear differential operators or, as we shall sometimes say, of adapted linear differential operator. Given two $C^m$ families $E$ and $E'$ of differentiable vector bundles on $\alpha\colon X\to S$, a linear differential operator carrying sections of $E$ to sections of $E'$ will be said to be a {\em $C^m$ family of linear differential operators} on $\alpha\colon X\to S$ if, when written in adapted coordinates and relative to adapted local trivializations of $E$ and $E'$, it involves only differentiation with respect to vertical coordinates and its coefficients are adapted functions. Thus an adapted linear differential operator carries adapted sections to adapted sections.

We will need a basic result, due to Kodaira and Spencer \cite{Kodaira_Spencer}, concerning the differentiability properties of solutions of differential equations
$Lu=v$, where $L$ is a $C^m$ family of differential operators and $v$ is an adapted section. Let $\alpha\colon X\to S$ be a $C^m$ family of compact differentiable manifolds, and let $E$ be a $C^m$ family of differentiable vector bundles on it. We denote by $\widetilde A (E)$ the vector space of adapted sections of $E$. Let $L\colon \widetilde A (E)\to\widetilde A (E)$ be a $C^m$ family of linear differential operators. A metric on $E$ will be said to be adapted if the inner product of any pair of adapted sections is an adapted functions. Adapted metrics always exist, and can for instance be constructed by gluing together flat local metrics by means of a partition of unity made up of adapted functions. Suppose an adapted metric is given on $E$, and one on the relative tangent bundle to $X\to S$. We denote by $\langle\ ,\ \rangle_s$ the inner product on $E_s$, and by $dV_s$ the volume form on $X_s$ coming from the metric. Consider the inner product
\be
(u,v)=\int_{X_s}\langle u,v\rangle_s\, dV_s\label{innprod}
\ee
on $A(E_s)$, the vector space of $C^\infty$ sections of $E_s$. We will say that $L$ is a {\em family of formally self-adjoint, strongly elliptic differential operators} if each $L_s$ is self-adjoint with respect to the inner product (\ref{innprod}), and strongly elliptic. Under these circumstances the kernel of $L_s$ is finite-dimensional, and there are linear operators
$$
F_s,G_s\colon A(E_s)\to A(E_s)\,,
$$
where $F_s$ is the orthogonal projection onto the kernel of $L_s$, and
\be
u=F_su+L_sG_su\,.\label{green}
\ee
for any $u\in A(E_s)$. We shall refer to $F_s$ and $G_s$, respectively, as the {\em harmonic projector} and {\em Green operator} associated to $L_s$ and to the chosen metrics.
The theorem of Kodaira and Spencer reads as follows.

\begin{theorem}[\cite{Kodaira_Spencer}, Theorem 5]\label{parameters}Let $m$ be a non-negative integer or $\infty$.
Let $\alpha\colon X\to S$ be a $C^m$ family of compact differentiable manifolds, and let $E$ be a $C^m$ family of differentiable vector bundles on $\alpha\colon X\to S$. Suppose $E$ and the relative tangent bundle to $X\to S$ are endowed with adapted metrics. Let $L\colon \widetilde A (E)\to\widetilde A (E)$ be a $C^m$ family of formally self-adjoint, strongly elliptic linear differential operators. Suppose the dimension of the kernel of $L_s$ is independent of $s$. Then the family of harmonic projectors $F=\{F_s\}_{s\in S}$ and the family of Green operators $G=\{G_s\}_{s\in S}$ are of class $C^m$, in the sense that $Fu$ and $Gu$ are adapted, for any adapted $u$.
\end{theorem}
To be precise, the statement proved by Kodaira and Spencer is slightly less general than the one we have given, in two respects. First of all, they deal only with families of the form $X_0\times S\to S$, where $X_0$ is a compact manifold. More importantly, they treat only the case $m=\infty$. These, however, are not serious difficulties. Since the statement of Theorem \ref{parameters} is local on $S$, the first is resolved by the following result.

\begin{lemma}\label{prodstruct}Let $\alpha\colon X\to S$ be a $C^m$ family of compact differentiable manifolds, let $E$ be a $C^m$ family of differentiable vector bundles on it, and let $s_0$ be a point of $S$. Then, if $U$ is a sufficiently small neighborhood of $s_0$, there is an isomorphism of $C^m$ families of differentiable manifolds between $\alpha^{-1}(U)\to U$ and the product family $X_{s_0}\times U\to U$. Moreover, if $q$ is the projection of $\alpha^{-1}(U)\cong X_{s_0}\times U$ to $X_{s_0}$, $E$ is isomorphic, as a $C^m$ family of differentiable vector bundles, to $q^*(E_{s_0})$.
\end{lemma}

This result is well known, at least for $m=\infty$. However, since the usual proof, which involves integrating lifts to $X$ of coordinate vector fields on $S$, breaks down for $m=0$, we shall sketch an alternate proof below. Before we do so, however, we notice that the second difficulty mentioned above is also non-existent, since the proof given by Kodaira and Spencer for their theorem, and in particular their crucial Proposition 1, work equally well, and virtually without changes, in our context, provided we substitute Lemma \ref{prodstruct} above, and in particular its second part, for their Lemma 1.

We now prove our Lemma \ref{prodstruct}. We shall use the following version ``with parameters'' of the inverse function theorem.

\begin{lemma}\label{inversepar}Let $U$ be an open subset of $\RR^n$, $V$ an open subset of $\RR^\ell$, $f\colon U\times V\to \RR^n$ a continuous function, and $m$ a non-negative integer or $\infty$. Write $x=(x_1,\dots,x_n)$ to indicate the standard coordinates in $U$ and $t=(t_1,\dots,t_\ell)$ to indicate the standard ones in $V$. Suppose that:
\begin{itemize}
\item[\sl i)]the function $f$ is $C^\infty$ in $x$ for any fixed $t$;
\item[\sl ii)]the function $f$ and all its derivatives, of any order, with respect to the $x$ variables, are $C^m$ functions of $x$ and $t$;
\item[\sl iii)]the Jacobian $\frac{\del f}{\del x}(x,t)$ is non-singular at a point $(x_0,t_0)\in U\times V$.
\end{itemize}
Set $F(x,t)=(f(x,t),t)$. Then there is an open neighborhood $A$ of $(x_0,t_0)$ such that $F(A)$ is open in $\RR^n\times \RR^\ell$ and $F$ induces a homeomorphism from $A$ to $F(A)$. Moreover, writing the inverse of this function under the form $(x,t)\mapsto (g(x,t),t)$, the function $g$ is $C^\infty$ in $x$, and $g$ and all its derivatives, of any order, with respect to the $x$ variables, are $C^m$ functions of $x$ and $t$.
\end{lemma}
The proof of this lemma is essentially the same as the one of the standard inverse function theorem, and will not be given here.

\begin{proof}[Proof of Lemma \ref{prodstruct}]The first step is to show that families of compact differentiable manifolds can be locally embedded in euclidean space. More precisely, we shall see that, if $U$ is a sufficiently small neighborhood of $s_0$, for sufficiently large $N$ there is a map $F\colon \alpha^{-1}(U)\to \RR^N$ such that the pair $({\bf 1}_U,F\times \alpha)$ is a morphism of $C^m$ families of differentiable manifolds from $\alpha^{-1}(U)\to U$ to $\RR^N\times U\to U$ which is a fiberwise embedding.
Denote by $B_r$ the ball of radius $r$ centered at the origin of $\RR^n$. Shrinking $S$, if necessary, we may find finitely many adapted charts $(\varphi_i,\alpha)\colon U_i\to B_2\times S$, where $U_i$ is an open subset of $X$, such that $X$ is covered by the open sets $(\varphi_i,\alpha)^{-1}(B_1\times S)$. Choose a non-negative $C^\infty$ function $\rho$ on $B_2$ which vanishes identically on the complement of $B_{4/3}$ and is identically equal to $1$ on $B_1$. Pulling this back via $\varphi_i$ we get a function on $U_i$, which we denote by $\rho_i$. Write $\varphi_i=(\varphi_{i,1},\dots,\varphi_{i,n})$, and denote by $\psi_{i,j}$ the function $\rho_i\varphi_{i,j}$, extended to zero on the complement of $U_i$ in $X$; clearly, $\psi_{i,j}$ is adapted. Denote by $\psi$ the map from $X$ to $\RR^M$ whose components are the $\psi_{i,j}$. The map $\psi$ is adapted, and its restriction to each fiber of $\alpha$ is a local embedding. In particular, there is an open neighborhood $W$ of the diagonal $\Delta$ in $X\times_S X$ such that $\psi(x)\neq\psi(y)$ for all $(x,y)\in W\smallsetminus \Delta$. Shrinking $S$ again, we may find finitely many adapted charts $(\xi_i,\alpha)\colon V_i\to B_2\times S$, where $V_i$ is an open subset of $X$, having the following property. We may select a set $I$ of pairs of indices such that $V_i\times_S V_j$ does not meet $\Delta$ if $(i,j)\in I$, and $X\times_S X\smallsetminus W$ is covered by the open sets $(\xi_i,\alpha)^{-1}(B_1\times S)\times_S (\xi_j,\alpha)^{-1}(B_1\times S)$ as $(i,j)$ varies in $I$. Denote by $\lambda_i$ the pullback of $\rho$ via $\xi_i$, extended to zero on the complement of $V_i$. By construction, if $(x,y)\in X\times_S X\smallsetminus W$, there is an index $i$ such that $\lambda_i(x)\neq 0$ but $\lambda_i(y)=0$. As a consequence, the $\lambda_i$, together with the $\psi_{i,j}$, are the components of an adapted map $X\to \RR^N$ which is a fiberwise embedding.

We may thus view $X$ as embedded in $\RR^N\times S$, and hence each fiber $X_s$ as embedded in $\RR^N$. There is a neighborhood $V$ of $X_{s_0}$ in $\RR^N$ which is diffeomorphic to a neighborhood of the zero section of the normal bundle to $X_{s_0}$ in $\RR^N$, and the projection to the zero section in the normal bundle yields a $C^\infty$ map $\eta\colon V\to X_{s_0}$. Clearly, $X_s\subset V$ when $s$ is close to $s_0$, and therefore, after suitably shrinking $S$, $\eta$ gives an adapted map $\beta\colon X\to X_{s_0}$. We wish to show that, possibly after further shrinking $S$, $(\beta,\alpha)\colon X\to X_{s_0}\times S$ is an isomorphism of $C^m$ families of differentiable manifolds. As the embedding of $X$ in $\RR^N$ is adapted, the tangent spaces to the $X_s$ vary continuously. Moreover, Lemma \ref{inversepar} implies that every point of $X_{s_0}$ has a neighborhood $U$ such that $\beta$ induces a diffeomorphism from $X_s\cap \beta^{-1}(U)$ to $U$ for $s$ near $s_0$. Let us see that, in fact, $\beta_{|X_s}$ is a diffeomorphism when $s$ is close to $s_0$. It suffices to prove injectivity. Suppose there are a sequence $\{s_n\}$ in $S$ converging to $s_0$ and sequences of points $x_n,y_n\in X_{s_n}$ such that $x_n\neq y_n$ but $\beta(x_n)=\beta(y_n)$; we may assume that $\{x_n\}$ and $\{y_n\}$ converge, respectively, to points $x,y\in X_{s_0}$. Letting $n$ go to infinity we find that $x=\beta(x)=\beta(y)=y$. This shows that, if $U$ is neighborhood of $x$ as above, $x_n$ and $y_n$ belong to $U$ for large $n$, and hence $\beta(x_n)\neq \beta(y_n)$, a contradiction. This proves the first part of Lemma \ref{prodstruct}.

The proof of the last statement is similar. We sketch it for $E$ a real vector bundle, the proof for a complex bundle being entirely analogous. First notice that, up to shrinking $S$, there exist enough adapted sections of $E$ to embed $E$ as a sub-vector bundle in a trivial bundle $\RR^K\times X_{s_0}\times S$, for some large $K$. The orthogonal projection with respect to the Euclidean metric of $\RR^K$ gives a surjective morphism of $C^\infty$ vector bundles on $X_{s_0}$ from $\RR^K\times X_{s_0}$ to $E_{s_0}$. As a consequence, we get a surjective morphism $\RR^K\times X_{s_0}\times S\to q^*(E_{s_0})$ which, when composed with the inclusion of $E\subset \RR^K\times X_{s_0}\times S$ yields a morphism of adapted vector bundles $E\to q^*(E_{s_0})$. The restriction of this to $\alpha^{-1}(U)$ is an isomorphism for any sufficiently small neighborhood $U$ of $s_0$.
\end{proof}

The notion of $C^m$ family of differentiable manifolds has a holomorphic counterpart in the one of {\em $C^m$ family of complex manifolds}. Formally, such a family is a surjective morphism $\alpha\colon X\to S$ of $C^m$ manifolds with the property that each fiber $X_s$ is a complex manifold, satisfying the following local triviality condition. For every $x\in X$ there is a $C^m$ diffeomorphism $\varphi:U\to V\times W$,
where $U$ is a neighborhood of $x$ in
$X$, $W$ is a neighborhood of $\alpha(x)$ in $S$, and $V$ is a ball centered at $0$ in some $\CC^h$, such that:
\begin{itemize}
\item[\sl i)]$\varphi$ is compatible with the projections to $S$;
\item[\sl ii)]$\varphi(x)=(0,\alpha(x))$;
\item[\sl iii)]$\varphi$ maps $U\cap \alpha^{-1}(s)$ biholomorphically onto $V\times \{s\}$ for every $s\in W$.
\end{itemize}
When all the fibers of $\alpha$ are curves we will speak of {\em $C^m$ family of curves}. A $C^m$ family of compact complex manifolds has a natural structure of $C^m$ family of differentiable manifolds, as a consequence of the following standard result.

\begin{lemma}Let $f(z_1,\dots,z_h,t_1,\dots,t_\ell)$ be a $C^m$
function of the complex variables $z_1,\dots,z_h$ and of the real variables
$t_1,\dots,t_\ell$ which is holomorphic in $z_1,\dots,z_h$. Then the partials of $f$, of any order,
with respect to the variables $z_1,\dots,z_h$, are $C^m$ functions of
$z_1,\dots,z_h,t_1,\dots,t_\ell$.
\end{lemma}
The lemma follows directly from Cauchy's integral formula
\bes
\dfrac{\del^{\sum k_i}}{\del z_1^{k_1}\dots\del z_h^{k_h}}
f(z_1,\dots,z_h,t_1,\dots,t_k)
=\frac{\prod k_i!}{(2\pi i)^h}
\int_{\genfrac{}{}{0pt}{}{\vert\zeta_i-z_i\vert=\varepsilon}{i=1,\dots,h}}\hskip-8pt
\frac{f(\zeta_1,\dots,\zeta_h,t_1,\dots,t_\ell)}{\prod(\zeta_i-z_i)^{k_i+1}}
d\zeta_1\dots d\zeta_h\,.
\ees
In fact, the right-hand side can be continuously differentiated $m$ times under the
integral sign with respect to the variables $t_i$.

\bigskip
A morphism of $C^m$ families of complex manifolds is just a morphism of $C^m$ families of differentiable manifolds which happens to be holomorphic on the fibers. If $\alpha\colon X\to S$ is a $C^m$ family of complex manifolds, we denote by $\widetilde{\mathcal O}_X$ the sheaf of $C^m$ functions on $X$ which are holomorphic along the fibers, and by $\widetilde{\Omega}^{p}_\alpha$ the $\widetilde{\mathcal O}_X$-module whose sections are the relative differential forms which restrict to a holomorphic $(p,0)$-form on each fiber. We also write $\widetilde{\mathcal A}^{p,q}_\alpha$ for the sheaf of adapted relative $(p,q)$-forms; it is a module over $\widetilde{\mathcal A}_X$, the sheaf of adapted functions. Now let $C$ be a curve. A {\em $C^m$ deformation} of $C$ is the datum of a pointed $C^m$ manifold $(S,s_0)$, a $C^m$ family of curves $\alpha\colon X\to S$, and an isomorphism $\varphi\colon C\overset\simeq\to X_{s_0}$. A {\em morphism of $C^m$ deformations} from $\alpha\colon X\to (S,s_0)$, $\varphi\colon C\to X_{s_0}$ to another deformation $\alpha'\colon X'\to (S',s'_0)$, $\varphi'\colon C\to X'_{s'_0}$ is a pair $(f,F)$ of $C^m$ morphisms fitting in a cartesian diagram
$$
\xymatrix{
X \ar[r]^F \ar[d]_\alpha & X' \ar[d]^{\alpha'}\\
(S,s_0) \ar[r]^f & (S',s'_0) \\
}
$$
such that $(f,F)$ is a morphism of $C^m$ families of curves and $F_{s_0}\circ\varphi=\varphi'$.
We shall show that the universal property of the Kuranishi family of $C$ holds also with respect to $C^m$ deformations.

\begin{proposition}\label{diffKura}Let $C$ be a curve of genus $g>1$, and let $\alpha\colon X\to (S,s_0)$, $\varphi\colon C\to X_{s_0}$, be a $C^m$ deformation of $C$. Let $\pi:\mathcal C\to(B,b_0)$, $\psi\colon C\to \mathcal C_{b_0}$, be a Kuranishi family for $C$. Then, for some open neighborhood $A$ of $s_0$, there is a morphism $(f,F)$ of $C^m$ deformations from $\alpha^{-1}(A)\to (A,s_0)$, $\varphi\colon C\to X_{s_0}$ to the Kuranishi family. This morphism is essentially unique, in the sense that any morphism of deformations from $\alpha^{-1}(A')\to (A',s_0)$, $\varphi\colon C\to X_{s_0}$, $A'$ a neighborhood of $s_0$, to the Kuranishi family, agrees with $(f,F)$ on $\alpha^{-1}(U)\to U$, where $U$ is a neighborhood of $s_0$.
\end{proposition}

\begin{proof}We claim that, possibly after shrinking $S$, we can find a cartesian diagram
\be
\xymatrix{
X \ar[r]^\Lambda \ar[d]_\alpha & \mathcal Y \ar[d]^{\xi}\\
(S,s_0) \ar[r]^\lambda & (H_{g,\nu},s'_0) \\
}\label{complexify}
\ee
where $(\lambda,\Lambda)$ is a morphism of $C^m$ families of complex manifolds and $\xi\colon\mathcal Y\to H_{g,\nu}$ is the universal family over the Hilbert scheme of $\nu$-canonically embedded genus $g$ curves, for some $\nu\ge 3$. The existence part of the proposition then follows from the standard universal property of the Kuranishi family, applied to $\mathcal Y\to H_{g,\nu}$.

To prove the existence of (\ref{complexify}) we proceed as follows. Let $u_{s_0}$ be a holomorphic $\nu$-fold differential on $X_{s_0}$, for some $\nu\ge 3$. We shall show that $u_{s_0}$ extends, possibly after shrinking $S$, to a section of $\mathcal L=(\widetilde{\Omega}^1_\alpha)^{\otimes \nu}$ on all of $X$. First of all, by Lemma \ref{prodstruct}, $u_{s_0}$ extends to a section $u$ of $(\widetilde{\mathcal A}^{1,0}_\alpha)^{\otimes \nu}=\mathcal L\otimes_{\widetilde{\mathcal O}_X}\widetilde{\mathcal A}_X$. Put adapted metrics on the $C^m$ families of differentiable vector bundles underlying $\mathcal L$ and $\widetilde{\mathcal A}^{1,0}_\alpha$. For each $s\in S$, we write $u_s$ for the restriction of $u$ to the fiber $X_s$, and $\delbar_s$ to indicate the $\delbar$ operator acting on sections of $\mathcal L\otimes_{\widetilde{\mathcal O}_X}\mathcal{A}_{X_s}$. We also let $\vartheta_s$ be the formal adjoint of $\delbar_s$, and $\Box_s=\delbar_s\vartheta_s$ the Laplace-Beltrami operator, both acting on $\mathcal L$-valued $(0,1)$-forms on $X_s$. The kernel of $\Box_s$ can be identified with $H^1(X_s,\mathcal O(K_{X_s}^\nu))$, where $K_{X_s}$ stands for the canonical bundle of $X_s$, and hence is reduced to zero, since $g>1$ and $\nu\ge 3$. In particular, the family of differential operators $\{\Box_s\}$ satisfies the assumptions of Theorem \ref{parameters}. It follows from (\ref{green}) that
$$
\delbar_su_s=\Box_sG_s\delbar_su_s=\delbar_s\vartheta_sG_s\delbar_su_s\,.
$$
By Theorem \ref{parameters}, $\{G_s\delbar_su_s\}$, and hence $v=\{v_s\}=\{\vartheta_sG_s\delbar_su_s\}$, are adapted. Since $\delbar_s(u_s-v_s)=0$ for any $s$, and $v_{s_0}=0$, we conclude that $u-v$ is a section of $\mathcal L$ extending $u_{s_0}$.

Now choose a basis for the space of holomorphic $\nu$-fold differentials on $X_{s_0}$ and, after suitably shrinking $S$, extend its elements to sections $u_0,\dots,u_N$ of $\mathcal L$ on $X$, where $N=(2\nu-1)(g-1)-1$. We may assume that the restrictions of these sections to $X_s$ constitute a basis for the space of holomorphic $\nu$-fold differentials on $X_s$, for any $s$. The sections $u_0,\dots,u_N$ then give a morphism of $C^m$ families of complex manifolds over $S$
$$
\xymatrix{
X \ar[rr] \ar[dr]_\alpha && \mathbb P^N\times S \ar[dl]^\eta\\
&S&\\
}
$$
such that $X_s\to\mathbb P^N$ is a $\nu$-canonical embedding for each $s$. To construct (\ref{complexify}), we take as $\lambda$ the map which associates to each $s\in S$ the point of $H_{g,\nu}$ corresponding to $X_s\hookrightarrow \mathbb P^N$. Recall that $\lambda(s)$ can be described as follows. The pullback homomorphism
\be
H^0(\mathbb P^N,\mathcal O(k))\to H^0(X_s,\mathcal O_{X_s}(k))\cong H^0(X_s,\mathcal O(K_{X_s}^{k\nu}))\label{hilbpoint}
\ee
is onto for every $k\ge 1$, by Noether's theorem. Its kernel is a point of the Grassmannian $\mathbb G$ of $\big(\binom{N+k}{N}-(2k\nu-1)(g-1)\big)$-planes in $\binom{N+k}{N}$-space. For large enough $k$, the Hilbert scheme $H_{g,\nu}$ is a subscheme of $\mathbb G$, and $\lambda(s)$ is just the kernel of (\ref{hilbpoint}). On the other hand, by construction, the homomorphism (\ref{hilbpoint}) is the fiber at $s$ of a homomorphism of $C^m$ vector bundles
$$
\eta_*\widetilde{\mathcal O}_{\mathbb P^N\times S}(k)\to \alpha_*\mathcal L^k\,.
$$
Since the kernel of this homomorphism is a $C^m$ vector subbundle of the trivial bundle on the left, the map $\lambda$ is $C^m$.

The map $\Lambda$ is easily constructed. If $x$ is a point of $X$, then
$$
\Lambda(x)=([u_0(x):\dots:u_N(x)],\lambda(s))\in X'\subset \mathbb P^N\times H_{g,\nu}\,.
$$
It is clear that $(\lambda,\Lambda)$ is a morphism of $C^m$ families.

It remains to prove the uniqueness property of $(f,F)$. Let $(f',F')$ be another morphism of deformations from $\alpha^{-1}(A')\to (A',s_0)$, $\varphi\colon C\to X_{s_0}$ to the Kuranishi family of $C$. We may suppose that $A'=A$. We may also suppose that $\operatorname{Aut}(C)$ acts on $\mathcal C$ and $B$, and that any isomorphism between fibers of $\pi$ comes from an automorphism of $C$. Assume that there is a sequence $\{x_n\}$ in $A$, converging to $s_0$, such that the restriction of $F$ to $X_{x_n}$ is different from the one of $F'$ for each $n$. We will show that this leads to a contradiction. In fact, what the assumption implies is that for each $n$ there is a non-trivial element $\gamma_n$ of $\operatorname{Aut}(C)$ such that $\gamma_nf(x_n)=f'(x_n)$ and $\gamma_nF$ coincides with $F'$ on $X_{x_n}$. As the automorphism group of $C$ is finite, we may suppose that all the $\gamma_n$ are equal to a fixed $\gamma\in\operatorname{Aut}(C)$. Then, passing to the limit for $n\to\infty$, we conclude that $\gamma F$ and $F'$ agree also on $X_{s_0}$. On the other hand, since we are dealing with morphisms of deformations, $F_{s_0}\circ\varphi=F'_{s_0}\circ\varphi$, which implies that $\gamma=1$, a contradiction.

\end{proof}
\np The notion of {\tei} structure carries over, with obvious changes, to the context of $C^m$ families of curves. An important corollary of Proposition \ref{diffKura} is then the following.
\begin{theorem}Let $\alpha\colon X\to S$ be a $C^m$ family of genus $g$ curves with {\tei} structure. Suppose $g\ge 2$. Let $\eta\colon \mathcal X_{g} \to T_g$ be the universal family on the genus $g$ {\tei} space. Let $f\colon S\to T_g$ be the map which associates to each point of $S$ the isomorphism class of the corresponding fiber, and let $F\colon X\to \mathcal X_g$ be the map whose restriction to the fiber $X_s$ is the unique isomorphism of curves with {\tei} structure from $X_s$ to $\eta^{-1}(f(s))$. Then the pair $(f,F)$ is a morphism of $C^m$ families of curves. In particular, $f$ is of class $C^m$.
\end{theorem}

\section{The theorem of Teichm\"uller}\label{teichthm}

In this section we discuss the following famous theorem of {\tei}.

\begin{theorem}\label{teich1}Let $g>1$.
The {\tei} space $T_g$ is homeomorphic to the unit ball in
$\CC^{3g-3}$.
\end{theorem}

\np
Our plan is to show how the proof of this theorem directly reduces to
the proof of two fundamental results, namely {\tei}'s uniqueness theorem
and the theorem of existence of solutions of the Beltrami equation.

Before we can proceed, we must introduce the notion of Beltrami differential. Let $S$ be a compact connected Riemann surface. We consider vector-valued differentials on $S$ which are locally of the form
$$
\mu=\nu\frac{\del}{\del z}\otimes d\overline{z},
$$
where $z$ is a local coordinate and $\nu$ is a measurable function, that is, measurable sections of $T_S\otimes \overline{K}_S$, where $T_S$ and $K_S$ are, respectively, the complex tangent and cotangent bundles to $S$. It makes sense to define a measurable function $|\mu|$ by setting it locally equal to $|\nu|$, since the latter is clearly independent of the choice of coordinate $z$.
A {\em Beltrami differential} on the Riemann surface $S$ is an $L^\infty$ section of $T_S\otimes \overline{K}_S$ whose norm $\Vert\mu\Vert=\sup_S|\mu|$ is strictly less than 1. We shall really need only a particular kind of Beltrami differentials, namely those which are $C^\infty$ everywhere, except at a finite number of points. These differentials will be called {\em admissible}.

To a Beltrami differential we associate a perturbed version of the $\delbar$ operator on $S$, by setting
$$
\delbar_\mu=\delbar-\mu\,,
$$
where $\mu$ acts on a function $f$ as
$$
\mu(f)=\nu\frac{\del f}{\del z}d\overline{z}\,.
$$
The corresponding {\em Beltrami equation} is
$$
\delbar_\mu u=0\,,
$$
that is,
\be
u_{\ov{z}}=\nu(z)u_z\,.\label{beltreq1}
\ee
The basic existence theorem for the Beltrami equation asserts that it has local solutions, in an appropriate generalized sense, which are homeomorphisms to open subsets of the complex plane. We will not need the full strength of this result, which is due to Morrey \cite{Morrey}, but just the fact that the same conclusion holds under the stronger hypothesis that $\mu$ is $C^\infty$. This is due to Korn \cite{Korn} and Lichtenstein \cite{Lichtenstein}, who more generally deal with the case when $\mu$ satisfies a H\"older condition; a simplified proof of their result was given by Bers \cite{Bers_7} and Chern \cite{Chern}, and can be found also in Chapter IV, section 8, of \cite{CourantHilbert}. Formally, the existence result we need is the following.

\begin{theorem}\label{localbeltrexistence}Let $\nu(z,t_1,\dots,t_n)$ be a $C^\infty$ function on
a neighborhood of the origin in $\CC\times\RR^n$. Suppose that
$|\nu(z,t)| <1$ for all values of $z$ and $t=(t_1,\dots,t_n)$. Then there exists a
$C^\infty$ function $w(z,t)$, also defined on a neighborhood of the origin, such that
$$
w_{\overline{z}}=\nu(z,t)w_z\,,\quad w_z(0,0)\neq 0\,.
$$
\end{theorem}			

We shall also need the following uniqueness result.

\begin{lemma}\label{localbeltruniqueness}Let $\nu(z)$ be a function on
a neighborhood of the origin in the complex plane which is $C^\infty$ except at a finite set $Z$. Suppose that $|\nu(z)|< 1$ for all $z$. Let $u$ be a homeomorphism from a neighborhood of the origin to an open subset of $\CC$ which solves the Beltrami equation (\ref{beltreq1}) away from $Z$. Let $f$ be a bounded function on a neighborhood of the origin which is once differentiable away from $Z$. Then $f$ satisfies (\ref{beltreq1}) away from $Z$ if and only if it is a holomorphic function of $u$.
\end{lemma}
The proof of the lemma is quite straightforward. Suppose first that $Z$ is empty, and let $w$ be the solution of the Beltrami equation provided by Theorem \ref{localbeltrexistence}. A simple chain rule computation shows that
$$
f_{\ov{z}}-\nu f_z=(1-|\nu |^2)\cdot f_{\ov{w}}\cdot \ov{w_z}\,.
$$
Thus $f$ is a holomorphic function of $w$ near the origin if and only if it is a solution of the Beltrami equation. This applies in particular to the function $u$. Moreover, since $u$ has an inverse, it is also the case that $w$ is a holomorphic function of $u$. This proves the lemma when $Z$ is empty. In the general case what the above argument shows is that $f$ solves the Beltrami equation away from $Z$ if and only if it is a holomorphic function of $u$ there. Holomorphicity of $f$ at points of $Z$ then follows from the Riemann extension theorem.

Theorem \ref{localbeltrexistence} in its parameterless version, that is, for $n=0$, and Lemma \ref{localbeltruniqueness} say, in particular, that a $C^\infty$ Beltrami differential $\mu$ on a Riemann surface $S$ defines on $S$ a new complex structure whose holomorphic functions are the solutions of the corresponding Beltrami equation.
The reader should be warned that $\delbar_\mu$ is {\em not} the $\delbar$ operator of this complex structure, but just its component of type $(0,1)$ (with respect to the original structure). A more suggestive, though cumbersome, notation for $\delbar_\mu$ could thus be $\delbar_\mu^{(0,1)}$.

Dependence on parameters in Theorem \ref{localbeltrexistence} is often not considered in the literature; a notable exception is \cite{Ahlfors_3}, to which we might refer for a proof. A cheap alternative is to appeal instead to Theorem \ref{parameters}. Here is how the argument goes. Since the problem is of a local nature, we may alter $\nu$ outside a neighborhood of the origin. Hence we may assume that $\nu$ is defined and $C^\infty$ on $\CC\times U$, where $U$ is a neighborhood of the origin in $\RR^n$, and that there is a positive $r$ such that $\nu$ vanishes for $|z|>r$. Thus we may view $\mu=\nu\,\del/\del z\otimes d\ov{z}$ as a family $\{\mu_t\}_{t\in U}$ of Beltrami differentials on $\PP^1$, vanishing outside the disk of radius $r$ centrered at $0$. Let $H$ be the hyperplane bundle on $\PP^1$; its smooth sections can be viewed as $C^\infty$ functions $f$ on $\CC$ such that $f/z$ extends in a $C^\infty$ way across $\infty$. Since $\delbar_{\mu_t}$ is the standard $\delbar$ operator in a neighborhood of $\infty$, for any smooth section $u$ of $H$, $\delbar_{\mu_t}u$ is a smooth $H$-valued $(0,1)$-form. Let $\vartheta_{\mu_t}$ be the formal adjoint of $\delbar_{\mu_t}$ with respect to (say) the Fubini-Study metric. Then the differential operator $L_t=\vartheta_{\mu_t}\delbar_{\mu_t}$ is self-adjoint and strongly elliptic. If $u$ is a section of $H$ such that $L_tu=0$, then $(\delbar_{\mu_t}u,\delbar_{\mu_t}u)=(u,L_tu)=0$, and hence $\delbar_{\mu_t}u=0$. In other words, the $L_t$-harmonic sections of $H$ are just those sections which are holomorphic with respect to the complex structure defined by $\mu_t$. The space of these sections has dimension $2$ for any $t$, by the Riemann-Roch theorem, so that Theorem \ref{parameters} applies. Pick a $\mu_0$-holomorphic section $u$ which vanishes simply at $0$, and set $w_t=F_tu$, where $F_t$ is the harmonic projector associated to $L_t$. Then $\delbar_{\mu_t}w_t=0$, and $w_t$ depends differentiably on $t$, by Theorem \ref{parameters}. In other words, $w(z,t)=w_t(z)$ has all the required properties.

\bigskip
Admissible Beltrami differentials originate, in particular, from the so-called admissible quasi-diffeomorphisms. An orientation-preserving homeomorphism $F\colon S\to S'$
between two compact Riemann surfaces which is a
diffeomorphism outside a finite set $Z\subset S$ is called a
{\it quasi-diffeomorphism}. Pick holomorphic coordinates
$z$ and $w$ around $p\in S\smallsetminus Z$ and $F(p)\in S'$,
respectively.
The condition that $F$ be an orientation-preserving diffeomorphism on
$S\smallsetminus Z$ tells us that, on $S\smallsetminus Z$, the Jacobian
determinant of $F$ is positive. Since, in local coordinates, the Jacobian is
$$
\vert w_z\vert ^2-\vert w_{\overline{z}}\vert ^2\,,
$$
the local function
$$
\nu(z)=\frac{w_{\overline{z}}}{w_z}
$$
is $C^\infty$ away from $Z$ and of absolute value less than 1.
It is a straightforward application of the chain rule to check that setting
$$
\mu_{_{F}}=\nu\frac{\del}{\del z}\otimes
d\overline z
$$
gives a well-defined section of $T_S\otimes\overline K_S$ which is $C^\infty$ away from $Z$. The quasi-diffeomorphism $F$ is said to be {\it admissible} if
$\Vert \mu_F\Vert <1$, i.e., if $\mu_F$ is a Beltrami differential. By Lemma \ref{localbeltruniqueness}, the
complex structure of $S'$ can be completely described in terms of the one of $S$ and of the differential $\mu_F$; a bounded function $f$ on an open subset of $S'$ is holomorphic if and only if $u=f\circ F$ is a solution of the Beltrami equation $\delbar_{\mu_F}u=0$ away from $Z$.
\begin{lemma}\label{mufmufinv}Let $F\colon S\to S'$ be a quasi-diffeomorphism. Then
$$
\Vert\mu_F\Vert=\Vert\mu_{F^{-1}}\Vert\,.
$$
If $F'\colon S\to S''$ is another quasi-diffeomorphism, then $\mu_{F'}=\mu_F$ if and only if $F'\circ F^{-1}\colon S'\to S''$ is holomorphic.
\end{lemma}
The second assertion of the lemma follows immediately from Lemma \ref{localbeltruniqueness}, while another elementary chain rule computation shows that $|\mu_F(p)|=|\mu_{F^-1}(F(p))|$ for any $p\in S$, thus proving the first assertion.

It is convenient to introduce the concept of \it dilatation
\rm for a quasi-diffeomorphism $F\colon S\to S'$. This is simply defined to be
$$
K[F]=\frac{1+\Vert \mu_{_{F}}\Vert }{1-\Vert \mu_{_{F}}\Vert }\,.
$$
It follows from Lemma \ref{mufmufinv} that a quasi-diffeomorphism and its inverse have
the same dilatation. It is also clear that
$F$ is admissible if and only if
$K[F]<\infty$.

\bigskip
We now turn to Theorem \ref{teich1}. Let $S$ be a reference Riemann surface of genus $g>1$. As we have announced, we shall rely on a fundamental result of {\tei}, the so-called {\tei} uniqueness theorem.
This provides a {\it canonical} representative for each isotopy class
of orientation-preserving
homeomorphisms $f\colon C\to S $, which has the following two remarkable
properties:
\begin{enumerate}
\item it is an admissible quasi-diffeomorphism;
\item away from the points where it fails to be smooth it can be locally described, {\it in a canonical way}, as a {\it real affine} transformation.
\end{enumerate}
Let $\omega$ be a holomorphic quadratic differential on $S$. If $\omega=f(z)dz^2$ is a local
expression for it in terms of a local coordinate $z$, we define the (singular)
volume form associated to $\omega$ to be
$$
dA_\omega=\frac{i}{2}|f|dz\wedge d\overline{z}
=|f|dx\wedge dy\,,
$$
where $x$ and $y$ stand for the real and imaginary parts of $z$. One immediately
checks that this definition is independent of the choice of local coordinate.
Now consider the
$(3g-3)$-dimensional space $H^0(S,\mathcal O(K^2_S))$ of holomorphic
quadratic differentials on $S$. A norm is introduced in this space by setting
$$
\Vert \omega\Vert =\int_S dA_\omega
$$
for any $\omega$ in $H^0(S,\mathcal O(K^2_S))$. Look at the unit
ball in $H^0(S,\mathcal O(K^2_S))$:
$$
B(K^2_S)=\{\omega\in
H^0(S,\mathcal O(K^2_S)):\Vert \omega\Vert <1\}\,.
$$
We are going to define a map
$$
\Phi\colon B(K^2_S)\to T_S=T_g
$$
and prove the following more precise version of (\ref{teich1}).

\begin{theorem}\label{teich2}$\Phi$ is a homeomorphism.
\end{theorem}

\np The construction of $\Phi$ is as follows. Pick $\omega\in B(K^2_S)$ and set $$
k=\Vert \omega\Vert <1\,.
$$
Let $Z$ be the set of zeroes of $\omega$. We introduce a new complex structure on $S\smallsetminus Z$. Near a point $p\in S\smallsetminus Z$ we may write
$$
\omega=(dz)^2\,,
$$
where $z$ is a holomorphic local coordinate, which is uniquely defined up to sign and the addition of a constant. We will say that $z$ is an {\it $\omega$-coordinate}. We define a new coordinate patch around $p$ by setting
\be
z'=\frac{z+k\overline{z}}{1-k}.
\label{teich.def}
\ee
A different choice of $z$ changes $z'$ by at most a sign
and the addition of a constant. Hence these coordinate patches define a new
holomorphic structure on $S\smallsetminus Z$.
We will show that this structure extends to one on all of $S$ by constructing
explicit coordinate patches around each point of $Z$; this extension
will clearly be unique. Let
$p$ be a zero of $\omega$, and $z$ a local coordinate around
$p$. Write $\omega=f(z)dz^2$ near $p$. We first treat the case when $f$ vanishes to even order $2n$ at $p$. By suitably changing coordinate, we may suppose that $\omega=d(z^{n+1})^2$. It follows that $z^{n+1}$ is an
$\omega$-coordinate at all points of $S$ sufficiently close to $p$, but different
from it, and hence, by virtue of (\ref{teich.def}), that
$$
\xi(z)=\frac{z^{n+1}+k\overline{z}^{n+1}}{1-k}
$$
is a local coordinate for the new complex structure on $S\smallsetminus Z$ at all
these points. Since $k<1$, the function
$1+k\frac{\overline{z}^{n+1}}{z^{n+1}}$
takes its values in the half-plane of complex numbers with positive real part,
where a single-valued determination of the $(n+1)$-st root function exists. Set
$$
\eta(z)
=z\left(\frac{1
+k\frac{\overline{z}^{n+1}}{z^{n+1}}}{1-k}\right)^{\frac{1}{n+1}}
$$
for $z\neq 0$ and $\eta(0)=0$. Since
$1+k\frac{\overline{z}^{n+1}}{z^{n+1}}$ is bounded, the function $\eta$ is continuous at $p$. We will now show that
$\eta$ is a homeomorphism between an open neighborhood of $p$ and an open
neighborhood of the origin in the complex plane, and we will take it as our new
coordinate around $p$. This is compatible with the new
complex structure on $S\smallsetminus Z$, since $\eta^{n+1}=\xi$. To prove our claim, by the ``invariance of domain'' theorem, it suffices to show that $\eta$ is injective. If
$\eta(z)=\eta(z')$, then $\xi(z)=\xi(z')$, and hence
$z^{n+1}={z'}^{n+1}$. It follows that
$$
1+k\frac{\overline{z}^{n+1}}{z^{n+1}}=1+k\frac{\overline{z'}^{n+1}}{{z'}^{n+1}}
$$
and therefore that $z=z'$, by the definition of $\eta$.
We now treat the case when $f$ vanishes to odd order $m$ at $p$. Let $q$ be the map
$\zeta\mapsto\zeta^2=z$. The pulled-back differential
$q^*\omega$ vanishes to order $2m+2$ at the origin, and by suitably changing coordinate we may suppose that
$q^*\omega=d(\zeta^{m+2})^2$. This means that $\zeta^{m+2}$ is an
$\omega$-coordinate at all points of $S$ close to $p$, but different from it. Thus
$$
\xi=\frac{\zeta^{m+2}+k\overline{\zeta}^{m+2}}{1-k}
$$
is a local coordinate for the new complex structure on $S\smallsetminus Z$ at
these points. We take as new coordinate at $p$ the function $\eta(z)$ whose
value is
$$
\eta(z)=z\left(\frac{1
+k\frac{|z|^{m+2}}{z^{m+2}}}{1-k}\right)^{\frac{2}{m+2}}
=\zeta^2 \left(\frac{1
+k\frac{\overline{\zeta}^{m+2}}{\zeta^{m+2}}}{1-k}\right)^{\frac{2}{m+2}}
$$
for $z\neq 0$ and zero for $z=0$. To prove that this extends the new complex
structure on $S\smallsetminus Z$ one proceeds as in the even
order case, after noticing that $\eta^{m+2}=\xi^2$. This finishes the
construction of the new complex structure of $S$.

We may now define the {\tei} map $\Phi\colon B(K^2_S)\to T_g$. We set
$$
\Phi(\omega)=[S_\omega, [f_\omega]]\,,
$$
where $S_\omega$ is the surface $S$, equipped with the new complex structure we have just described, and $f_\omega\colon S_\omega\to S$ is the
set-theoretic identity. The homeomorphism $f_\omega$
is called the {\tei} map associated to $\omega$, and is an admissible quasi-diffeomorphism. In fact, by (\ref{teich.def}), $\nu_{f_\omega}=-k$, and hence
$$
|\mu_{f_\omega}|=k=\Vert\omega\Vert<1
$$
everywhere.

It is crucial to observe that $S_\omega$ and $f_\omega$
depend continuously on $\omega\in B(K^2_S)$. More precisely, if
we denote by $\mathcal S$ the disjoint union of all the $S_\omega$ and by $f$ the
map from $\mathcal S$ to $S\times B(K^2_S)$ which sends $x\in
S_\omega$ to $(f_\omega(x),\omega)$, what we just did was to put on
$\mathcal
S\to B(K^2_S)$ a structure of
continuous family of compact genus $g$ Riemann surfaces for which
$f$ is a topological trivialization. The map $f$ thus endows the family
$\mathcal S \to B(K^2_S)$ with a {\tei} structure. By the universal property of {\tei} space, discussed in the preceding section, the
map $\Phi$ is therefore {\it continuous}. We see that this rather simple proof of the continuity of $\Phi$ is a direct consequence of our definition of {\tei} space in terms of Kuranishi families.

The uniqueness theorem of {\tei} asserts that, among all admissible
quasi-diffeomor\-phisms isotopic to it, $f_\omega\colon S\to
S$ is one with
{\it minimal dilatation}, and that it is uniquely characterized by
this property.
\begin{theorem}{\bf (Teichm\"uller's uniqueness theorem). \rm}Let
$S$ be a genus $g$ Riemann surface with
$g>1$, and let $\omega$ be a holomorphic quadratic differential on $S$. Assume that $\Vert \omega\Vert =k<1$. Let
$f\colon S_\omega\to S $ be an admissible
quasi-diffeomorphism which is isotopic to $f_\omega$ (i.e., isotopic to the
identity). Then
$$
K[f]\geq K[f_\omega]= \frac{1+k}{1-k}\,,
$$
and equality holds if and only if $f=f_\omega$.
\end{theorem}

\np
We will not prove the theorem here, but we refer to one of the many proofs which exist in the literature (see, for example, \cite{Abikoff}, \cite{Gardiner_1}, or \cite{Lehto_1} and the bibliography therein).

Assuming Teichm\"uller's uniqueness theorem, we shall first
prove that the map $\Phi\colon B(K^2_S)\to T_g$ is {\it injective}. Suppose that
$\Phi(\omega_1)=\Phi(\omega_2)$; in other words, that there are an isomorphism
$\varphi\colon S_{\omega_1}\to S_{\omega_2}$ and
an isotopy $f_{\omega_1}\sim f_{\omega_2}\circ \varphi$. Set $k_i=
\Vert \omega_i\Vert $, $i=1,2$.
The first part of Teichm\"uller's uniqueness
theorem, together with Lemma \ref{mufmufinv}, tells us that
$$
\frac{1+k_2}{1-k_2}=K[ f_{\omega_2}]=K[f_{\omega_2}\circ\varphi]\geq
K[ f_{\omega_1}]=\frac{1+k_1}{1-k_1}\,.
$$
Reversing the roles of $\omega_1$ and $\omega_2$, we obtain $k_1=k_2$, and the second
part of Teichm\"uller's uniqueness theorem implies that
$$
f_{\omega_1}=f_{\omega_2}\circ \varphi\,.
$$
Let $z_i$ be an $\omega_i$-coordinate on $S$, for $i=1,2$.
By Lemma \ref{mufmufinv} we have
$$
k_1\frac{\del}{\del z_1}\otimes d\overline{z}_1
=\mu_{f_{\omega_1}^{-1}}
=\mu_{f_{\omega_2}^{-1}}
=k_2\frac{\del}{\del z_2}\otimes d\overline{z}_2\,.
$$
But $k_1=k_2$, so that
$$
\overline{\frac{\del z_2}{\del z_1}}
=\frac{\del z_2}{\del z_1}\,.
$$
Since $z_2$ is a holomorphic function of $z_1$, this implies that
$z_2=c\cdot z_1+b$, where $b$ and $c$ are constants, and $c$ is real.
Hence $\omega_2=(dz_2)^2=c^2\cdot (dz_1)^2=c^2\cdot\omega_1$.
As $\Vert \omega_1\Vert =k_1=k_2=\Vert \omega_2\Vert $, we obtain that $c^2=1$,
which proves the injectivity of $\Phi$.

\bigskip
We are now going to conclude the proof
of Theorem \ref{teich2}, and hence of Theorem
\ref{teich1}. The first step is the following.

\begin{proposition}\label{phiclosed}The {\tei} map $\Phi$ is closed.
\end{proposition}
\begin{proof}Let $\{\omega_n\}$ be a sequence in
$B(K^2_S)$. Suppose the sequence $y_n=\Phi(\omega_n)$
converges to $y\in T_g$. What must be proved is that a subsequence of
$\{\omega_n\}$ converges in $B(K^2_S)$.
Set $y=[C, [f]]$, where $f$ is a diffeomorphism. As a
neighborhood of $y$ we take the basis $B$ of a Kuranishi family
$\pi\colon \mathcal C\to B$
for $C$.
We can assume that there is a
$C^\infty$ trivialization $(F,\pi)\colon \mathcal C\to S\times B$ such that
$F_{y}=f$. We may also assume that $\{y_n\}\subset B$.
Set $f_n=F_{y_n}$, $\mu=\mu_{_f}$ and
$\mu_n=\mu_{_{f_n}}$. Since $f$ and $f_n$ are $C^\infty$, the
Beltrami differentials $\mu$ and $\mu_n$ are also $C^\infty$; since $F$ is
$C^\infty$, $\{\mu_n\}$ converges uniformly to $\mu$.
Teichm\"uller's uniqueness theorem gives
$$
\frac{1+\Vert \mu\Vert }{1-\Vert \mu\Vert }=\underset{n\to
\infty}{\text{lim}}\left(\frac{1+\Vert \mu_n\Vert }{1-\Vert \mu_n\Vert }\right)=
\underset{n\to \infty}{\text{lim}}K[f_{n}]\geq \underset{n\to
\infty}{\text{lim}}K[f_{\omega_n}]
=\underset{n\to
\infty}{\text{lim}}\left(\frac{1+\Vert \omega_n\Vert }{1-\Vert \omega_n\Vert }\right).
$$
Therefore, for any constant $c$ with $\Vert \mu\Vert <c<1$, one has that
$\Vert \omega_n\Vert <c$, if $n$ is large enough. Thus a subsequence of
$\{\omega_n\}$ converges.
\end{proof}
\np
The last ingredient of the proof of Theorem \ref{teich2} is the following.
\begin{proposition}\label{teichconnect}The {\tei} space $T_g$ is connected.
\end{proposition}
\np
We can immediately see that this implies Theorem \ref{teich2}.
In fact, Proposition \ref{phiclosed}, in addition to showing that $\Phi(B(K^2_S))$ is closed in $T_g$, also shows that $\Phi$ gives a homeomorphism between
$B(K^2_S)$ and $\Phi(B(K^2_S))$. But
then the ``invariance of domain'' theorem says that
$\Phi(B(K^2_S))$ is open in $T_g$, since
$B(K^2_S)$ and $T_g$ are differentiable manifolds of
dimension $6g-6$.

We now turn to Proposition \ref{teichconnect}. An immediate consequence of Theorem \ref{localbeltrexistence} and Lemma \ref{localbeltruniqueness} is the following.
\begin{lemma}\label{beltrapar}Let $\mu_t$ be a family of smooth Beltrami
differentials on $S$, where $t$ varies in an interval $I\subset \RR$. Suppose
$\mu_t$ depends smoothly on $t$, in the sense that it is $C^\infty$ on $S\times I$. Then there are a differentiable family $\xi\colon Y\to I$ of
Riemann surfaces and a differentiable trivialization $S\times I\to Y$ such
that the Beltrami differential $\mu_{F_t}$ associated to $F_t$ is $\mu_t$.
\end{lemma}
\np
To prove Proposition \ref{teichconnect}, denote by $x_0$ the base point
$[S,[1]]$ of $T_S=T_g$, let $x=[C,[f]]$ be another element of $T_g$, where
$f\colon C\to S$ is a diffeomorphism, and set $\mu_t=t\mu_{f^{-1}}$. By Lemma \ref{beltrapar}, there is a differentiable family of curves with {\tei} structure over an
interval $I$ whose fiber at $t$ is $[Y_t,[F_t^{-1}]]$. This comes from a differentiable map $\gamma\colon I\to T_g$, by the universal
property of {\tei} space. To prove connectedness, we just have to show that $\gamma(0)=x_0$ and $\gamma(1)=x$; in other words, that $[Y_0,[F_0^{-1}]]=[S,[1]]$ and $[Y_1,[F_1^{-1}]]=[C,[f]]$. Since, by construction, $\mu_{F_0^{-1}}=0$ and $\mu_{F_1^{-1}}=\mu_{f^{-1}}$, this follows from the second part of Lemma \ref{mufmufinv}. This concludes the proof of Proposition \ref{teichconnect}, and therefore of {\tei}'s Theorem.

\bibliographystyle{amsplain}

\providecommand{\bysame}{\leavevmode\hbox to3em{\hrulefill}\thinspace}

\end{document}